\documentstyle{amsppt}

\NoBlackBoxes
\NoRunningHeads

%\textwidth=6.0in
%\textheight=10.0in

%\font\ref=cmr9
%\font\refit=cmti9
%\font\refbf=cmbx9

\nopagenumbers

\input amstex
\magnification=\magstep1

\def\ell{{\text{ell}}}

\def\h1{\hat{\bold 1}}

\def\Ua{U_q(\tilde\g)}
\def\U2{{\Ua}_2}
\def\g{\frak g}

\def\Z{\Bbb Z}

\def\d{\partial}

\def\<{\langle}
\def\>{\rangle}

\def\e{\varepsilon}

\topmatter
\title Factorization of differential operators, 
quasideterminants, and nonabelian Toda field equations
\endtitle

\author Pavel Etingof, Israel Gelfand, and Vladimir Retakh
\endauthor

\address
\newline
I.~G. :  Department of Mathematics, Rutgers University, New Brunswick,
NJ 08903
\newline
P.~E. and V.~R.: Department of Mathematics, Harvard University, Cambridge, MA
02138
\endaddress

\email
\newline
P.~E. : etingof\@ math.harvard.edu
\newline
I.~G. : igelfand\@ math.rutgers.edu
\newline
V.~R. : retakh\@ math.harvard.edu
\endemail

\endtopmatter

\centerline{\bf Abstract}

We integrate nonabelian Toda field equations \cite{Kr} 
for root systems of types $A$, $B$, $C$, 
for functions with values in any associative algebra. The solution
is expressed via quasideterminants introduced in \cite{GR1},\cite{GR2},
\cite{GR4}. In the appendix we review some results concerning noncommutative 
versions of other classical integrable equations. 

\centerline{\bf Introduction}
 
Nonabelian Toda equations for the root system $A_{n-1}$ 
were introduced by Polyakov (see \cite{Kr}). 
They are equations with respect to $n$ unknowns 
$\phi=(\phi_1,...,\phi_n)\in A[[u,v]]$, where $A$ is some associative
(not necessarily commutative) algebra with unit:
$$
\frac{\d}{\d u}\left(\frac{\d\phi_j}{\d v}\phi_j^{-1}\right)=
\cases\phi_2\phi_1^{-1},&j=1\\
\phi_{j+1}\phi_j^{-1}-\phi_j\phi_{j-1}^{-1},&2\le j\le n-1\\
-\phi_n\phi_{n-1}^{-1},&j=n\endcases
\tag 1
$$

Suppose that $A$ is a $*$-algebra, i.e. it is equipped with 
an involutive antiautomorphism $*:A\to A$. Then, setting in (1) $\phi_{n+1-i}=
(\phi_i^*)^{-1}$, we obtain a new system of equations. If $n=2k$,
 we get the nonabelian Toda system for root system
$C_k$:
$$
\frac{\d}{\d u}\left(\frac{\d\phi_j}{\d v}\phi_j^{-1}\right)=
\cases\phi_2\phi_1^{-1},&j=1\\
\phi_{j+1}\phi_j^{-1}-\phi_j\phi_{j-1}^{-1},&2\le j\le k-1\\
(\phi_k^*)^{-1}\phi_k^{-1}-\phi_k\phi_{k-1}^{-1},&j=k\endcases
\tag 2
$$
If $n=2k+1$, we get the nonabelian Toda system for root system $B_k$:
$$
\frac{\d}{\d u}\left(\frac{\d\phi_j}{\d v}\phi_j^{-1}\right)=
\cases\phi_2\phi_1^{-1},&j=1\\
\phi_{j+1}\phi_j^{-1}-\phi_j\phi_{j-1}^{-1},&2\le j\le k\\
(\phi_k^*)^{-1}\phi_{k+1}^*-\phi_{k+1}\phi_k^{-1},&j=k+1,\endcases
\tag 3
$$
where $\phi_{k+1}^*=\phi_{k+1}^{-1}$. 

Quasideterminants were introduced in \cite{GR1}, as follows. 
Let $X$ be an $m\times m$-matrix over $A$. 
For any $1\le i,j\le m$, let 
$r_i(X)$, $c_j(X)$ be the i-th row and the j-th column of $X$. 
Let $X^{ij}$ be the submatrix of $X$ obtained by removing 
the i-th row and the j-th column from $X$. For a row 
vector $r$ let $r^{(j)}$ be $r$ without the j-th entry.
For a column vector $c$ let $c^{(i)}$ be $c$ without the i-th entry. 
Assume that $X^{ij}$ is invertible. Then the quasideterminant 
$|X|_{ij}\in A$ is defined by the formula
$$
|X|_{ij}=x_{ij}-r_i(X)^{(j)}(X^{ij})^{-1}c_j(X)^{(i)},\tag 4
$$
where $x_{ij}$ is the $ij$-th entry of $X$.
 
In this paper we will use quasideterminants to integrate 
nonabelian Toda equations for root systems of types $A,B,C$. We do not know yet
how to generalize these results to root systems of types D-G, and to
affine root systems. 

Our method of integration, which is based on interpreting the Toda flow
as a flow on the space 
of factorizations of a fixed ordinary differential operator.
This method and explicit solutions of 
Toda equations are well known in the commutative case 
(see \cite{LS}; see also \cite{FF} and references therein;
for the one variable case see \cite{Ko}).  
 
\centerline{\bf 1. Factorization of differential operators} 

Let $R$ be an associative 
algebra over a field $k$ of characteristic zero,
and $D:R\to R$ be a $k$-linear derivation. 

Let $f_1,...,f_m$ be elements of $R$. By definition, the Wronski
matrix $W(f_1,...,f_m)$ is 
$$
W(f_1,...,f_m)=\left(\matrix f_1&...&f_m\\
Df_1&...&Df_m\\ ..&...&..\\ D^{m-1}f_1&...&D^{m-1}f_m\endmatrix\right)
$$
We call 
a set of
elements $f_1,...,f_m\in R$ nondegenerate if $W(f_1,...,f_m)$ is invertible. 

Denote by $R[D]$ the space of
polynomials of the form $a_0D^n+a_1D^{n-1}+...+a_n$, $a_i\in R$.  
It is clear that any element of $R[D]$ defines a linear operator on $R$. 

{\bf Example:} $R=C^\infty(\Bbb R),D=\frac{d}{dt}$. In this case, $R[D]$ 
is the algebra of differential operators on the line. 

By analogy with this example, 
we will call elements of $R[D]$ differential operators.

We will consider operators of the form $L=D^n+a_1D^{n-1}+...+a_n$.
We will call such $L$ an operator of order $n$ with highest coefficient 1. 
Denote the space of all such operators by $R_n(D)$. 

\proclaim{Theorem 1.1}
(i) Let $f_1,...,f_n\in R$ be a nondegenerate
set of elements. Then there exists a unique 
differential operator $L\in R[D]$ 
of order $n$ with highest coefficient $1$, such that $Lf_i=0$
for $i=1,...,n$. It is
given by the formula
$$
Lf=|W(f_1,...,f_n,f)|_{n+1,n+1}.\tag 1.1
$$

(ii) Let $L$ be of order $n$ with highest coefficient 1, and
$f_1,...,f_n$ be a set of solutions of the equation $Lf=0$, 
such that for any $m\le n$ the set
of elements $f_1,...,f_m$ is nondegenerate. 
Then $L$ admits a factorization $L=(D-b_n)...(D-b_1)$, where
$$
b_i=(DW_i)W_i^{-1},\ W_i=|W(f_1,...,f_i)|_{ii}.\tag 1.2
$$
\endproclaim

\demo{Proof} (i) We look for $L$ in the form
$L=D^n+a_1D^{n-1}+...+a_n$. From the equations $Lf_i=0$ 
it follows that 
$$
(a_n,...,a_1)=-(D^nf_1,...,D^nf_n)W(f_1,...,f_n)^{-1}.
$$
By definition,
$$
\gather
|W(f_1,...,f_n,f)|_{n+1,n+1}=
D^nf-(D^nf_1,...,D^nf_n)W(f_1,...,f_n)^{-1}(f,Df,...,D^{n-1}f)^T=\\
D^nf+(a_n,...,a_1)(f,Df,...,D^{n-1}f)^T=Lf.
\endgather
$$

(ii) We will prove the statement by induction in $n$. 
For $n=1$, the statement is obvious. Suppose it is valid for the differential 
operator $L_{n-1}$ of order $n-1$ with highest coefficient 1, 
which annihilates $f_1,...,f_{n-1}$ (by (i), it exists and is unique). 
Set $b_n=(DW_n)W_n^{-1}$, and consider the operator 
$\tilde L=(D-b_n)L_{n-1}$. It is obvious that 
$\tilde Lf_i=0$ for $i=1,...,n-1$. Also, by (i) 
$$
\tilde Lf_n=(D-b_n)L_{n-1}f_n=(D-b_n)W_n=0.
$$
Therefore, by (i), $\tilde L=L$. 
$\square$\enddemo

Now consider the special case: $R=A[[t]]$, where $A$ is an associative  
 algebra over $k$, 
and $D=\frac{d}{dt}$ (here $t$ commutes with everything).
In this case, it is easy to show that nondegenerate sets of
solutions of $Lf=0$ exist, and are in 1-1 correspondence with 
elements of the group $GL_n(A)$, via $\bold f=(f_1,...,f_n)\to W(\bold f)(0)$. 

It is clear that two different sets of solutions of the 
equation $Lf=0$ can define the same factorization of $L$. However, 
to each factorization $\gamma$ of $L$ we can assign a set 
$\bold f_\gamma=(f_1,...,f_n)$ of
solutions of $Lf=0$, which gives back the factorization $\gamma$ 
under the correspondence of Theorem 1.1(ii). This 
set is uniquely defined by the condition 
that the matrix $W(\bold f_\gamma)(0)$ is lower triangular 
with 1-s on the diagonal. 

Here is a formula for computing $\bold f_\gamma$, which is well known
in the commutative case.  

\proclaim{Proposition 1.2}
If $\gamma$ has the form 
$$
L=(D-(Dg_n)g_n^{-1})...(D-(Dg_1)g_1^{-1}),
$$
where $g_i(0)=1$, then $\bold f_\gamma=(f_1,...,f_n)$, where
$$
f_j(t)=\int_{0}^t\int_0^{t_1}...\int_0^{t_{j-2}}
g_1(t)g_1(t_1)^{-1}g_2(t_1)g_2(t_2)^{-1}...g_j(t_{j-1})dt_{j-1}...dt_2dt_1,
\tag 1.3
$$
where $\int_0^u(\sum a_it^i)dt:=\sum a_i\frac{u^{i+1}}{i+1}$. 
\endproclaim

\demo{Proof} It is easy to see that if $\bold f=(f_1,...,f_n)$,
with $f_j$ given by (1.3), then $W(\bold f)(0)$ is 
strictly lower triangular. So it remains to show that 
$f_j$ is a solution of the equation $L_jf=0$, 
where $L_j=(D-(Dg_j)g_j^{-1})...(D-(Dg_1)g_1^{-1})$. 

We prove this by induction in $j$. The base of induction is clear, 
since from (1.3) we get $f_1=g_1$. Let us perform the induction step. 
By the induction assumption, From (1.3), we have
$$
f_j(t)=g_1(t)\int_0^tg_1(s)^{-1}h(s)ds,
$$
where $h$ obeys the equation 
$(D-(Dg_j)g_j^{-1})...(D-(Dg_2)g_2^{-1})h=0$.
Thus, we get
$$
(D-(Dg_1)g_1^{-1})f_j=h.
$$
This proves that $L_jf_j=0$.
$\square$\enddemo

Now consider an application of these results to the noncommutative 
Vieta theorem \cite{GR3}. Let $A$ be an associative  algebra.
Call a set of elements $x_1,...,x_n\in A$ generic if 
their Vandermonde matrix $V(x_1,...,x_n)$ ($V_{ij}:=x_j^{i-1})$
is invertible. 

Consider an algebraic equation
$$
x^n+a_1x^{n-1}+...+a_n=0.\tag 1.4
$$
with $a_i\in A$. Let $x_1,...,x_n\in A$ be solutions 
of (1.4) such that $x_1,...,x_i$ form a generic set
for each $i$. Let $V(i)=V(x_1,...,x_i)$, and 
$y_i=|V(i)|_{ii}x_i|V(i)|_{ii}^{-1}$.  

\proclaim{Theorem 1.3}\cite{GR3} (Noncommutative Vieta theorem)
$$
a_r=(-1)^r\sum_{i_1<...<i_r}y_{i_r}...y_{i_1}.\tag 1.5
$$
\endproclaim

\demo{Proof} (Using differential equations.) Consider the differential 
operator with constant coefficients in $R=A[[t]]$:
$$
L=D^n+a_1D^{n-1}+...+a_n.
$$
We have solutions $f_i=e^{tx_i}$ of the equation $Lf=0$, and 
for any $i$ the set $f_1,...,f_i$ is nondegenerate, since its Wronski matrix
$W(i)$ is of the form 
$$
W(i)=V(i)\text{diag}(e^{tx_1},...,e^{tx_i}).\tag 1.6
$$
Thus, by Theorem 1(ii), the operator
$L$ admits a factorization
$$L=(D-b_n)...(D-b_1),$$ where $b_i=(D|W(i)|_{ii})|W(i)|_{ii}^{-1}$. 
Substituting (1.6) into this equation, we get
$$
b_i=|V(i)|_{ii}x_i|V(i)|_{ii}^{-1}=y_i.\tag 1.7
$$
Since $[D,b_i]=0$, we obtain the theorem. 
\enddemo

\vskip .1in 

\centerline{\bf 2. Toda equations}

In this section we will solve equations (1),(2),(3).

Let $M=\{(g,h)\in A[[t]]\oplus A[[t]]: g(0)=h(0)\in A^*\}$,
where $A^*$ is the set of invertible elements of $A$.
We have the following obvious proposition, which 
says that solutions of Toda equations are uniquely determined by initial 
conditions.  

\proclaim{Proposition 2.1} The assignment 
$\phi(u,v)\to (\phi(u,0),\phi(0,v))$ is 
a bijection between the set of solutions 
of the Toda equations (1) and the set $M$.
\endproclaim

Let $B$ be an algebra over $k$, 
$R=B[[v]]$, $D=\frac{d}{d v}$. 
For $\phi=(\phi_1,...,\phi_n)$, where $\phi_i\in B[[v]]$ are invertible, 
define $L^\phi_i\in R[D]$ by
$$
L_i^\phi=(D-(D\phi_i)\phi_i^{-1})...(D-(D\phi_1)\phi_1^{-1}).\tag 2.1
$$

Now set $B=A[[u]]$. Then $B[[v]]=A[[u,v]]$, so 
for any $\phi=(\phi_1,...,\phi_n)$, with all 
$\phi_i\in A[[u,v]]$ invertible, we can define $L_i^\phi$. 

\proclaim{Proposition 2.2} 
A vector-function $\phi(u,v)$ is a solution
of Toda equations (1) if and only if
$$
\frac{\d L_i^\phi}{\d u}=-\phi_{i+1}\phi_i^{-1}L_{i-1}^\phi, i\le n-1;
\ \frac{\d L_n^\phi}{\d u}=0.\tag 2.2
$$ 
\endproclaim

\demo{Proof} Let $\phi$ be a solution 
of the Toda equations. Set $b_i=(D\phi_i)\phi_i^{-1}$. We have 
$L_i^\phi=L_i=(D-b_i)...(D-b_1)$. 
Therefore, for $i\le n-1$ we have 
$$
\gather
\frac{\d L_i}{\d u}=-\sum_{j=1}^{i}
(D-b_i)...(D-b_{j+1})\phi_{j+1}\phi_j^{-1}
(D-b_{j-1})...(D-b_1)\\
+\sum_{j=1}^{i-1}
(D-b_i)...(D-b_{j+2})\phi_{j+1}\phi_{j}^{-1}
(D-b_{j})...(D-b_1).\tag 2.3
\endgather
$$
But 
$$
(D-b_{j+1})\circ \phi_{j+1}\phi_j^{-1}=
\phi_{j+1}\phi_j^{-1}(D-b_j).
$$
Therefore, the right hand side of (2.3)
equals $-\phi_{i+1}\phi_i^{-1}L_{i-1}$. 
The same argument shows that for $i=n$ the derivative $\frac{\d L_n}{\d u}$ 
vanishes.

Conversely, it is easy to see that equations (2.2) imply (1). 
$\square$\enddemo

Now we will compute the solutions of Toda equations explicitly. 
Let $\eta_1,...,\eta_n,\psi_1,...,\psi_n\in A[[t]]$ 
be such that $\eta_i(0)=\psi_i(0)\in A^*$. 
We will find the solution 
of the following initial value problem for Toda equations:
$$
\phi_i(u,0)=\psi_i(u),\phi_i(0,v)=\eta_i(v)\tag 2.4
$$
Proposition 2.1 states that this solution exists and is unique.

Let $g_i(v)=\eta_i(v)\eta_i(0)^{-1}$. Let $\bold f=
(f_1,...,f_n)$ be given by formula (1.3) in terms of $g_i$.
Define the lower triangular matrix $\Delta(u)$ whose 
entries are given by the formula
$$
\gather
\Delta_{ij}(u)=\int_0^u\int_0^{t_1}...\int_0^{t_{i-j-1}}
\psi_i(t_{i-j})\psi_{i-1}^{-1}(t_{i-j})\psi_{i-1}(t_{i-j-1})...
\psi_j^{-1}(t_1)\psi_j(u)dt_{i-j}...dt_1\endgather
$$
Let 
$\bold f^u=(f_1^u,...,f_n^u)$ be defined by the formula
$$
\bold f^u=\bold f\Delta(u).
$$
Then we have

\proclaim{Theorem 2.3} The solution of the problem (2.3) is given by the 
formula
$$
\phi_i(u,v)=|W(f_1^u,...,f_i^u)|_{ii}\tag 2.5
$$
\endproclaim

\demo{Proof} Let $L_i^u=L_i^{\phi(u,v)}$, where $\phi(u,v)$ is defined
by (2.5). 
By Theorem 1.1(ii), we have $L_i^uf_j^u=0$ for $j\le i$.
Differentiating this equation with respect to $u$, we get 
$$
L_i^u\frac{\d f_j^u}{\d u}+\frac{\d L_i^u}{\d u}f_j^u=0.\tag 2.6
$$
On the other hand, it is easy to see that
$$
\Delta'(u)=\Delta(u)\Theta(u), 
$$
where 
$$
\Theta_{ij}(u)=\cases \psi_i^{-1}(u)\psi_i'(u)&i=j\\
1&i=j+1\\ 0\text{ otherwise }\endcases
$$
This implies that 
$$
\frac{\d f_i^u}{\d u}=f_{i+1}^u+f_i^u\psi_i^{-1}\psi_i'.
$$
Therefore, $L_i^u\frac{\d f_j^u}{\d u}=0$ for $i\ge j+1$, 
and 
$$
L_i^u\frac{\d f_j^u}{\d u}=L_i^uf_{i+1}^u=|W(f_1^u,...,f_{i+1}^u)|_{i+1i+1}=
\phi_{i+1}. 
$$
Thus, from (2.6) we get
$$
\frac{\d L_i^u}{\d u}f_j^u=-\phi_{i+1}\phi_i^{-1}L_{i-1}f_j^u,\ i\ge j.
$$
By Theorem 1.1(i), this implies 
equations (2.2), which are equivalent to Toda equations. 

It is obvious that $\phi(u,v)$ satisfies the required initial conditions. 
The theorem is proved.
$\square$\enddemo

If initial conditions (2.4) satisfy the symmetry 
property $\phi_{n+1-i}=(\phi_i^*)^{-1}$, then we obtain
a solution of the initial value problem for equations (2) (for even $n$), and
(3) (for odd $n$). This gives a complete description of solutions 
of systems of equations (1),(2),(3). 

In the case $\phi(u,0)=1$, 
the Toda flow can be interpreted as a flow on the space of 
factorizations of a differential operator.
Namely, let $L$ be a differential operator of order 
$n$ with highest coefficient 1. Let $F(L)$ be the space of factorizations of 
$L$. Let $N^-_n$ be the group of strictly lower triangular 
matrices over $A$. It is easy to see that the map $\pi: F(L)\to N_n^-$,
given by $\pi(\gamma)=W(\bold f_\gamma)(0)$, is a bijection.
We will identify $F(L)$ with $N_n^-$ using $\pi$. 

Let $\gamma(u)=\gamma(0)e^{uJ_n}$ be a curve on $N_n^-$
generated by the 1-parameter subgroup $e^{uJ_n}$, where
$J_n$ is the lower triangular nilpotent Jordan matrix
($(J_n)_{ij}=\delta_{i,j+1}$).
Let $L=(D-b_n^u)...(D-b_1^u)$ be the factorization of $L$ corresponding 
to the point $\gamma(u)$. Let $\phi_i(u,v)$ 
be such that $(D\phi_i)\phi_i^{-1}=b_i^u$,
and $\phi_i(u,0)=1$. (here $D=\frac{\d}{\d v}$).
Then we have the following Corollary from Theorem 2.3.

\proclaim{Corollary 2.4} $\phi=(\phi_1,...,\phi_n)$ is a solution 
of Toda equations (1), and all solutions with
$\phi(u,0)=1$ are obtained in this way. 
\endproclaim

\demo{Proof}
For the proof it is enough to observe that if $\psi_i(u)=1$ then
$\Theta(u)=J_n$, and $\Delta(u)=e^{uJ_n}$. 
$\square$\enddemo

{\bf Remark.} An analogous statement can be made for equations (2) and (3). 
In this case, instead of an arbitrary differential operator $L\in R_n(D)$
one should consider a selfadjoint (respectively, skew-adjoint) 
operators, instead of the group $N_n^-$ -- 
a maximal nilpotent subgroup in 
$Sp_n(A)$ (respectively, $O_n(A)$), and 
instead  of $J_n$ the sum of simple root elements in the Lie algebra of this 
group. In the commutative case, such a picture of the Toda flow is known 
for all simple Lie groups \cite{FF}. 

Consider now the statement of Theorem 2.3 for $i=1$. In this case we have 
$$
\phi_1(u,v)=f_1^u(v)=\sum_{i=1}^nf_i^0(v)\Delta_{i1}(u)
\tag 2.7
$$

Thus, we get 

\proclaim{Corollary 2.4} \cite{RS} 
If $\phi$ is a solution of the Toda equations
then $\phi_1(u,v)=\sum_{i=1}^n p_i(u)q_i(v)$, where
$p_i,q_i$ are some formal series. 
\endproclaim

Now let us discuss infinite Toda equations. These are equations (1) with 
$n=\infty$. These equations allow to express $\phi_i$ recursively 
in terms of $f=\phi_1$, which can be done using quasideterminants:
$$
\phi_i=|Y_i(f)|_{ii}, \ Y_i(f):=
W\biggl(f,\frac{\d f}{\d u},...,\frac{\d^{i-1}f}{\d u^{i-1}}
\biggr),\tag 2.8
$$
where $\frac{\d^if}{\d u^i}$ are regarded as functions of $v$
for a fixed $u$. 
This formula is easily proved by induction. It appears in \cite{GR2},
section 4.5 in the case $u=v$, and in \cite{RS} in the 2-variable case. 

Formula (2.8) can be used to give another expression for the general 
solution of finite Toda equations, which appears in \cite{RS}. Indeed,
we have the following easy proposition.

\proclaim{Proposition 2.5} Let $f\in A[[u,v]]$ be such that 
$Y_1(f),...,Y_n(f)$
are invertible matrices. 
Then $|Y_{n+1}(f)|_{n+1n+1}=0$ if and only if 
$f$ is ``kernel of rank n'', i.e.
$$
f=\sum_{i=1}^n p_i(u)q_i(v).\tag 2.9
$$
\endproclaim

Thus, taking $f=\phi_1$ of the form (2.9), with $Y_1,...,Y_n$ invertible, 
and using formula (2.6)
we will get a solution $(\phi_1,...,\phi_n)$ of the finite Toda system of 
length $n$. It is not difficult to show that in this way one gets 
all possible solutions. 

{\bf Example.} Consider the nonabelian Liouville equation
$$
\frac{\d }{\d u}\biggl(\frac{\d \phi}{\d v}\phi^{-1}\biggr)=
(\phi^*)^{-1}\phi^{-1},\tag 2.10
$$
which is a special case of (2) for $k=1$. 
Consider the initial value problem
$$
\phi(0,v)=\eta(v),\phi(u,0)=\psi(u), \eta(0)=\psi(0)=a.\tag 2.11
$$
By Theorem 2.3, we get the following formula for the solution:
$$
\phi(u,v)=\eta(v)\biggl(a^{-1}+
\int_0^v\int_0^u\eta^{-1}(t)(\eta^*)^{-1}(t)a^*(\psi^*)^{-1}(s)\psi^{-1}(s)
dsdt\biggr)\psi(u).\tag 2.12
$$
For example, if $\psi(u)=1$, we get 
$$
\phi(u,v)=\eta(v)\biggl(1+u\int_0^v\eta^{-1}(t)(\eta^*)^{-1}(t)dt\biggr).\tag 2.13
$$
In the commutative case, these formulas coincide with the standard formulas
for solutions of the Liouville equation.

\vskip .1in

\centerline{\bf Appendix: Noncommutative soliton equations} 

In this appendix we will review some results about noncommutative versions 
of classical soliton hierarchies (KdV,KP). These results are mostly known 
or can be obtained by a trivial generalization of the corresponding 
commutative results, but they have never been exposed systematically. 

We will follow Dickey's book \cite{D}.

{\bf A1. Nonabelian KP hierarchy.}

Nonabelian KP hierarchy is defined in the same way as 
the usual KP hierarchy \cite{D}. 
Let 
$$
L=\d+w_0\d^{-1}+w_1\d^{-2}+....\tag A1
$$ 
be a formal pseudodifferential operator.
Here $\d=\frac{d}{dx}$, $w_i\in A[[x]]$, where $A$ is an associative
(not necessarily commutative) algebra 
with $1$. Consider the following infinite system of differential
equations:
$$
\frac{\d L}{\d t_m}=[B_m,L], B_m=(L^m)_+,\tag A2
$$
where for a pseudodifferential operator $M$, we denote by $M_+,M_-$  
the differential and the integral parts of $M$ 
(i.e. all terms with nonnegative, respectively negative,  
powers of $\d$). For brevity we will write $L^m_\pm$ instead of $(L^m)_\pm$.

Each of the differential equations (A2) defines a formal flow 
on the space $P$ of pseudodifferential operators of the form (A1). 
Indeed, $[L^m_+,L]=-[L_-^m,L]$, and the order of $[L_-^m,L]$ is at most zero. 

\proclaim{Proposition A1} The flows defined by equations 
(A2) commute with each other. 
\endproclaim

\demo{Proof} The proof is the same as the proof of Proposition (5.2.3) in 
\cite{D}. Namely, one needs to check the zero curvature condition
$$
\frac{\d B_m}{\d t_n}-\frac{\d B_n}{\d t_m}-[B_m,B_n]=0,\tag A3
$$
which is done by a direct calculation given in \cite{D}. 
\enddemo

The hierarchy of flows defined by (A2) for $m=1,2,3,...$ 
is called the noncommutative KP hierarchy.

{\bf A2. Noncommutative KdV and nKdV hierarchies.}

As in the commutative case, we can restrict the KP hierarchy to 
the subspace $P_n\subset P$ of all operators $L$ such that $L^n$ is 
a differential operator, i.e. $L_-^n=0$ \cite{D}. 
This subspace is invariant under the KP flows, 
as $\frac{\d L_-^n}{\d t_m}=[L^m_+,L^n_-]_-$. 
The space $P_n$ can be identified with the space $D_n$ of
differential operators of the form
$$
M=\d^n+u_2\d^{n-2}+...+u_n,\tag A4
$$
by the map $L\to M=L^n$. The KP hierarchy induces a hierarchy of flows
on $D_n$ called the nKdV hierarchy:
$$
\frac{\d M}{\d t_m}=[M^{m/n}_+,M].\tag A5
$$
Among these flows, the flows corresponding to $m=nl$, $l\in\Z_+$, are
trivial, but the other flows are nontrivial. 

For $n=2$ the nKdV hierarchy is the usual KdV hierarchy. 
The first two nontrivial equations of this hierarchy are
$$
u_{t_1}=u_x,\ u_{t_3}=\frac{1}{4}(u_{xxx}+3u_xu+3uu_x).\tag A6
$$
The second equation is the noncommutative KdV equation. 

{\bf A3. Finite-zone solutions of the nKdV hierarchy.}
 
Denote the vector fields of the noncommutative nKdV hierarchy 
by $V_m$, $m\ne nl$. Let $V=\sum_{i=1}^p a_iV_i$, $a_p\ne 0$,
 be a finite linear 
combination of $V_i$ with coefficients $a_i\in l$.
Let $S(a)$, $a=(a_1,...,a_p)$ be the set of stationary points 
of $V$ in $D_n=A[[x]]^{n-1}$. 

The set $S(a)$ can be identified with $A^{(n-1)p}$. 
Indeed, $S(a)$ is the subset of $D_n$ defined by the differential equations 
$V(M)=0$, which have the form  
$$
\frac{d^pu_i}{d x^p}=F_i(u_j,u_j',...,u_j^{(p-1)}).\tag A7
$$
So $u_i$ are uniquely determined by their first $p$ derivatives at $0$.

Since the vector field $V$ commutes with the nKdV hierarchy, $S(a)$ is 
invariant under the flows of this hierarchy. Using the identification
$S(a)\to A^{(n-1)p}$, we can rewrite each of the nKdV flows as a
system of ordinary differential equations on $A^{(n-1)p}$. 
Thus, solutions of nKdV equations belonging to $S(a)$ can be 
computed by solving ordinary differential equations. Such solutions 
are called finite-zone solutions. 

{\bf Remark.} 
If $M\in S(a)$ then $[Q(M^{1/n})_+,M]=0$, where $Q(x)=\sum a_ix^i$. 
Thus, we have two commuting differential operators $M$ and $N=Q(M^{1/n})_+$
in 1 variable. If the algebra $A$ is finite dimensional over its center
then there exists a nonzero polynomial $R(x,y)$ with coefficients in the 
center of $A$, such that $R(M,N)=0$. This polynomial defines an algebraic 
curve, called the spectral curve. The operator $M$ can then be computed 
explicitly using the method of Krichever \cite{Kr2}. 
Similarly, all nKdV equations resticted to the space of such operators can be 
solved explicitly in quadratures. However, if
$A$ is infinite-dimensional over its center, the polynomial $R$ does 
not necessarily exist, and we do not know any way of computing $M$ explicitly. 
 
{\bf A4. Multisoliton solutions of the noncommutative KP hierarchy.}

Here we will construct N-soliton solutions of the KP and KdV hierarchies
in the noncommutative case. We will use the dressing method. 
In the exposition we will closely follow Dickey's book \cite{D}.

We will now contruct a solution of the 
KP hierarchy. Let $t=(t_1,t_2,...)$
Consider the formal series
$$
\xi(x,t,\alpha)=(x+t_1)\alpha+t_2\alpha^2+...
+t_r\alpha^r+...,\alpha\in A.\tag A8
$$
Fix $\alpha_1,...,\alpha_N,\beta_1,...,\beta_N, a_1,...,a_N\in A$, and
set
$$
y_s(x,t)=e^{\xi(x,t,\alpha_s)}+a_s
e^{\xi(x,t,\beta_s)}.
\tag A9
$$
Define the differential operator of order $N$, with highest coefficient $1$, 
by  
$$
\Phi f(x)=|W(y_1,...,y_N,f)|_{N+1,N+1}(x),\tag A10
$$
where the derivatives in the Wronski matrix are taken with respect to $x$,
and we assume that the functions $y_1,...,y_N$ are a generic set
(in the sense of Chapter 1). 
Set 
$$
L=\Phi\d\Phi^{-1}.\tag A11
$$

\proclaim{Proposition A2} The operator-valued series
$L(x,t)$ is a solution 
of the KP hierarchy.
\endproclaim

\demo{Proof} The proof is the same as the proof of Proposition 5.3.6 in 
\cite{D}.  
\enddemo

Such solutions are called N-soliton solutions. 

Proposition A2 can be used to construct N-soliton solutions 
 of the nKdV hierarchy.
Namely, we should restrict the above construction to the case when 
$\beta_k=\e_k \alpha_k$, where $\e_k$ is an 
$n$-th root of unity. In this case it 
can be shown as in \cite{D} that the operator $M=L^n=\Phi\d^n\Phi^{-1}$
is a differential operator of order $n$, and it is a solution 
 of the nKdV hierarchy. 

As an example, let us consider the N-soliton solutions of the KdV 
hierarchy.  
In this case, we may set $t_{2k}=0$, and we have
$$
y_s=e^{\xi(x,t,\alpha_s)}+a_se^{-\xi(x,t,\alpha_s)}.\tag A12
$$

\proclaim{Proposition A3} Let 
$$
b_i=(\d W_i)W_i^{-1}, \ W_i:=|W(y_1,...,y_i)|_{ii}.\tag A13
$$
Then the function
$$
u(x,t)=2\d (\sum_{i=1}^Nb_i).\tag A14a
$$
is a solution of the noncommutative KdV hierarchy.
\endproclaim

\demo{Proof} Let $\Phi=\d^N+v_1\d^{N-1}+...$. 
Then from the equation $(\d^2+u)\Phi=\Phi\d^2$ we obtain 
$u=-2\d v_1$, and from Theorem 1.1(ii) $v_1=-\sum_i b_i$, where
$b_i$ are given by (A13). This implies (A14a). 
\enddemo

Another formula for $u$, which is equivalent to (A14a), is
the following. Let $Y(y_1,...,y_N)$ be the matrix which coincides with the 
Wronski matrix $W(y_1,...,y_N)$ 
except at the last row, where it has $y_i^{(N)}$ 
instead of $y_i^{(N-1)}$. Let $Y_N=|Y(y_1,...,y_N)|_{NN}$. 
Then the function $u$ given by (A14a) can be written as
$$
u=2\d(Y_NW_N^{-1}),\tag A14b.
$$

In the commutative case formulas (A14a),(A14b) 
reduce to the classical formula
$$
u=2\d^2\ln[\text{det}W(y_1,...,y_N)].\tag A15
$$ 

In particular, we can obtain N-soliton solutions of the 
noncommutative KdV equation
$u_t=\frac{1}{4}(u_{xxx}+3u_xu+3uu_x)$. For this purpose 
set $t_i=0$, $i\ne 3$, and $t_3=t$. 
Then we get 
$$
y_s=e^{\alpha_s x+\alpha_s^3 t}+a_se^{-\alpha_s x-\alpha_s^3 t},\tag A16
$$
and the solution $u(x,t)$ is given by (A14a),(A14b).

For example, consider the 1-soliton solution. According to (A14a), it has the 
form
$$
u=2\frac{\d}{\d x}
[(e^{\alpha x+\alpha^3 t}-ae^{-\alpha x-\alpha^3 t})\alpha
(e^{\alpha x+\alpha^3 t}+ae^{-\alpha x-\alpha^3 t})^{-1}].\tag A17
$$
In the commutative case, it reduces to the well known solution
$$
u=\frac{2\alpha^2}{\text{cosh}^2(\alpha x+\alpha^3 t-c)}, \ c=\frac{1}{2}\ln a,
\tag A18
$$
-- the solution corresponding to the solitary wave
which was observed by J.S.Russell in August of 1834.

\enddocument